\documentclass[12pt,twoside]{article}

\usepackage{amsmath,amsfonts,amsthm,amssymb,amsopn,amstext,amscd}
\usepackage{enumerate}
\usepackage{graphicx}
\usepackage{float}
\usepackage{setspace} %For the space between lines
\usepackage{enumitem}
\usepackage{hyperref}
\allowdisplaybreaks

\newcommand{\N}{\mathbb{N}}
\newcommand{\R}{\mathbb{R}}

\newcommand{\eps}{\varepsilon}
\newcommand{\cuad}{\begin{flushright}\vspace{-2ex}$\Box$\vspace{-2ex}\end{flushright}}

\newenvironment{Prf}[1][\unskip]{%
\par
\noindent
{\large\textbf{Proof of #1}}\newline
\vspace{-2ex}\noindent{}\newline}\cuad

\newtheorem{thm}{Theorem}
\newtheorem{prop}[thm]{Proposition}
\newtheorem{cor}[thm]{Corollary}
\newtheorem{lem}[thm]{Lemma}
\newtheorem{remark}{Remark}
\newtheorem{ex}{Examples}
\begin{document}
\pagenumbering{arabic}
%\pagestyle{myheadings}
%\singlespace
%\onehalfspace

%%%%%%%%%%%%%%%%%%%%%%%%%%%%%%%%%%%%%%%%%%%%%%%%%%%%%%%%%%%%%%%%%%%%%
\title{Branching processes in varying environment with generation-dependent immigration\footnote{This is the preprint version of the following paper published in the journal \emph{Stochastic Models}
(see the official journal website at \url{https://doi.org/10.1080/15326349.2019.1575754}):\\
M. Gonz\'alez, G. Kersting, C. Minuesa, I. del Puerto. Branching processes in varying environment
with generation-dependent immigration. Stochastic Models, 35(2): 148-166, 2019. DOI:
10.1080/15326349.2019.1575754}}
\author{M. Gonz\'alez, G. Kersting, C. Minuesa, I. del Puerto}

\date{July 23, 2018}

\maketitle
%%%%%%%%%%%%%%%%%%%%%%%%%%%%%%%%%%%%%%%%%%%%%%%%%%%%%%%%%%%%%%%%%%%%%

%%%%%%%%%%%%%%%%%%%%%%%%%%% ABSTRACT %%%%%%%%%%%%%%%%%%%%%%%%%%%%%%%%
\section*{Abstract}
A branching process in varying environment with generation-dependent immigration is a modification of the standard branching process in which immigration is allowed and the reproduction and immigration laws may vary over the  generations. This flexibility makes the process more appropriate to model real populations due to the fact that the stability in the reproductive capacity and in the immigration laws are not usually fulfilled. In this setting, we study the extinction problem, providing a necessary and sufficient condition for the certain extinction of the population. The asymptotic behaviour of the model is analysed for those processes with critical offspring distributions, according with the classification established in \cite{kersting-2017}, and when the immigration means stabilize to a positive value. More specifically, we establish that the asymptotic distribution of the process -under a suitable normalization- belongs to the gamma distribution family.

%%%%%%%%%%%%%%%%%%%%%%%%%%%%%%%%%%%%%%%%%%%%%%%%%%%%%%%%%%%%%%%%%%%%%

\vspace*{2ex}

\noindent\textbf{Keywords}: branching process; varying environment;  inhomogeneous immigration; extinction; asymptotic distribution.

%%%%%%%%%%%%%%%%%%%%%%% INTRODUCTION %%%%%%%%%%%%%%%%%%%%%%%%%%%%%%%%
\section{Introduction}
%%%%%%%%%%%%%%%%%%%%%%%%%%%%%%%%%%%%%%%%%%%%%%%%%%%%%%%%%%%%%%%%%%%%%

Since the appearance of the pioneer branching model, the Bienaym\'e-Galton-Watson process (BGWP), and motivated by the complexity of the problems arisen in diverse fields, new branching processes have been proposed with the aim of giving an answer to them. In particular, special attention has been paid to the problem of modelling migratory movements. Within the class of discrete time models, it is worthwhile to mention the branching process with immigration, the branching process with time-dependent immigration, the branching process with immigration only after empty generation, branching processes with migration or the branching process with immigration in varying or random environment (for their description, see the pioneer works \cite{Heathcote-1965,Heathcote-1966}, \cite{Foster-Williamson-1971}, \cite{Pakes-1971}, \cite{Nagaev-Khan-1980} and \cite{Key-1987}, respectively, or \cite{CBPs} for a review of them) and the age-dependent branching process with immigration or its generalization, the Sevast'yanov branching process with immigration (see \cite{jagers-68} and \cite{Sebastyanov-1957}, respectively, for their description), in continuous time among others.

There are several works where some of those models were considered for modelling in real situations. For instance, in \cite{Grishechkin-1992} a multitype continuous-time branching process with immigration is used to describe a retrial queue, where the birth and death of the individuals represent the arrival and departure of a customer, respectively, and the situations when a customer is immediately served after the arrival into the system are modelled through the immigration. In the context of cell kinetics, a two-type age-dependent branching process with emigration is used to model the emigration of oligodendrocyte cells cultured in vitro out of the field of observation in  \cite{Yakovlev-Stoimenova-Yanev-2008}; more recently, in \cite{art-Dposterior}, a controlled two-type branching process with binomial control is used to describe its discrete branching structure  as a result of an embedding of the aforesaid model. Another example is the application of an age-dependent branching process with time-inhomogeneous immigration for modelling the progression of Leukemia in mice in \cite{Hyrien-Yanev-Jordan-2015}; indeed, this process models the number of Leukemia cells in the blood and the immigration process represents the influx of cells from other tissues of the body.

This work is concerned with the branching process in a varying environment with generation-dependent immigration (BPVEI). This model represents a twofold generalization of the BGWP: on the one hand, we allow the reproductive capacity of the individuals to change over the generations and on the other hand, an inhomogeneous immigration is included into the probability model. These features make the process more appealing from a practical standpoint. Firstly, it is more natural to assume that the probability distribution governing the reproduction, referred as offspring distribution, may change over the time because the reproduction capacity of the individuals might be affected by factors such as, for instance, the resource supply or the weather conditions, which may vary over seasons. An example of such situations is the genus \emph{Bicyclus} of butterflies of the sub-Saharan Africa, whose reproduction capacity changes from the wet season to the dry-season due to the adaptation to the environment (see \cite{Haccou}, p.49, for further details). Secondly, the interest of considering an inhomogeneous immigration is due to in practical situations it is unknown whether the distribution of the influx of individuals in the population is constant or changes over the generations, and hence, the assumption of an homogeneous immigration might lead to a misspecification of the model. An example of such populations is a renewal process for cell populations of the central nervous system which comprise stem cells, progenitor cells and terminally differentiated cells. In such populations, only the stem and progenitor cells have reproductive capacity, whereas the terminally differentiated cells arise as a result of the maturation of progenitor cells. Since any of these types of cells could die at any time, an immigration of stem cells -and their development into progenitor cells- supplies the amount of cells needed to maintain the tissues. Moreover, the reproduction and immigration rates of stem cells is believed to depend on the renewing tissue and on the time according to the experimental evidences, and hence the need of an inhomogeneous immigration in the model (see \cite{Yakovlev-Yanev-2006} for further details).

While a great number of papers on branching processes in varying environment (BPVE), 
%which evolves as the BGWP with the only difference of that the offspring distribution might depend on the generation, 
or on branching processes with inhomogenous immigration is available (see, for instance, \cite{Fearn-1972}, \cite{Souza-1995},  \cite{kersting-2017} or \cite{BPVEs} for the first model and \cite{Foster-Williamson-1971} and \cite{Rahimov-1987, Rahimov-2007, Rahimov-2017,Rahimov-2018} for the second one), in relation to the research on the BPVEI one only can find the works \cite{Gao-Zhang-2015} and \cite{ispany-2016}. The former concerns a central limit theorem and a law of the iterated logarithm for the BPVEI in the supercritical case according to the classification in \cite{kersting-2017}, whereas in the second paper the author obtained a Feller diffusion approximation when assuming certain convergence of the parameters of both the offspring and immigration laws in the critical case. The paper \cite{Mitov-Omey-2014} also deals with a BPVE, but with homogeneous immigration and considering a very specific offspring distributions for such a study.

The main goals of this work are to analyse the extinction problem and the asymptotic behaviour for the general family of BPVEIs. Previously to the study of those questions, we establish some basic properties on the moments and the probability generating functions involved in the model. In relation to the extinction, it is important to highlight that although 0 is not an absorbing state of the process, it is possible to state conditions for the population size to reach the state 0 and to remain there forever. As we illustrate, this result also holds true when either the offspring distribution or the immigration laws are homogeneous. In order to examine the asymptotic behaviour of the process, we focus our attention on the family of BPVEIs with critical offspring distributions according to the classification in \cite{kersting-2017}. For those processes, we prove their convergence in distribution -once suitably normalized- to a gamma distribution when the immigration means and the normalized second factorial moments of the offspring distributions converge to positive values.

Apart from this introduction, this paper consists of three sections and an appendix. In Section \ref{sec:model}, we describe the probability model and some basic properties regarding the moments and probability generating functions of the process. Section \ref{sec:extinction} is devoted to the extinction problem and Section \ref{sec:lim-behaviour} is dedicated to the study of the asymptotic behaviour of the process. In order to ease the reading, we dedicate an appendix to the proofs of the results.

%%%%%%%%%%%%%%%%%%%%%%% MODEL  %%%%%%%%%%%%%%%%%%%%%%%%%%%%%%%%%%%%%%
\section{Probability model and basic properties}\label{sec:model}
%%%%%%%%%%%%%%%%%%%%%%%%%%%%%%%%%%%%%%%%%%%%%%%%%%%%%%%%%%%%%%%%%%%%%

In order to define a BPVEI, let us consider two independent families of $\N_0$-valued random variables, $\{X_{nj}:\ n\in\N_0;\ j\in\N\}$ and $\{I_{n}:\ n\in\N_0\}$. Assume that both of them are made up by independent random variables such that for each $n\in\N_0$ fixed, $X_{nj}$, $j\in\N$, are identically distributed according to the probability generating function (p.g.f.) $f_n(s)=\sum_{k=0}^\infty f_n[k]s^k$, with $f_n(0)<1$, and $I_n$ is distributed according to the p.g.f. $h_n(s)=\sum_{k=0}^\infty h_n[k]s^k$, with $h_n(0)<1$, for $n\in\N_0$. The process $\{Z_n\}_{n\in\N_0}$ defined recursively as
\begin{equation*}\label{def:model}
Z_0=0,\quad Z_{n+1}=\sum_{j=1}^{Z_{n}+I_n}X_{nj},\quad n\in\N_0,
\end{equation*}
is called \emph{branching process  in a varying environment  $v=\{f_0,f_1,f_2,\ldots\}$ with generation-dependent immigration} and initial value $Z_0=0$. The distributions $f_n=\{f_n[k]\}_{k\in\N_0}$ and $h_n=\{h_n[k]\}_{k\in\N_0}$ are called offspring distribution or reproduction law and immigration law, respectively, of the $n$-th generation. Notice that we use the same notation for the probability distributions and the corresponding p.g.f.s.

Intuitively, $Z_n$ represents the number of individuals in the $n$-th generation, $I_n$ is the number of immigrants in the $n$-th generation and $X_{nj}$ the number of offspring of the $j$-th progenitor in the $n$-th generation. Observe that the assumption $Z_0=0$ indicates that the population starts with the offspring of the immigrants in the initial generation. Moreover, in the development of this model one distinguishes two phases in each generation: the immigration phase, when new individuals migrate into the population, and the reproduction phase, when the progenitors have their offspring. Notice that we define the process as a controlled branching process in varying environment (see \cite{LMJ}) with an inhomogeneous control function $\phi_n(k)=k+I_n$, with $k,n\in\N_0$.
%The previous conditions mean that the immigration phases are characterized by the fact that they are independent of each others and of all the reproduction phases and the probability distribution governing the migration mechanism may depend on the generation. In each generation, the number of progenitors that participate in the reproduction phase is defined by the immigrants together with the individuals at the previous generation. Moreover, each of these progenitors reproduces  independently of the others and according to the same probability distribution, which is known as offspring distribution and which might be different from the offspring distribution at other generations. 

It is easy to verify that a BPVEI is an inhomogeneous Markov chain. We shall fix the notation for the remainder of the paper. Let us denote by $F_n$ the p.g.f. of the variable $Z_n$, $n\in\N_0$, and $f_{k,n}=f_{k+1}\circ\ldots\circ f_n$, for $k=-1,0,\ldots, n-1$, with the convention that $f_{n,n}(s)=s$, $s\in [0,1]$. Let us write $m_n=E[X_{n1}]$ and $\sigma_n^2=Var[X_{n1}]$ for the offspring mean and variance in generation $n$, respectively, and assume that both of them are positive and finite; analogously, denote by $\alpha_n=E[I_{n}]$ and $\beta_n^2=Var[I_{n}]$ the expectation and the variance of the number of immigrants in the $n$-th generation, which are assumed to satisfy $0<\alpha_n<\infty$ and $0<\beta_n^2<\infty$. For $n\in\N_0$, we also write $\mu_n=\prod_{i=0}^n m_i$, and $\nu_n=f_n''(1)/f_n'(1)^2$, with the convention $\mu_{-1}=1$.

%%%%%%%%%%%%%%%%%%%% P.G.F. AND MOMENTS %%%%%%%%%%%%%%%%%%%%%%%%%%%%%%
%\section{Probability generating functions and moments}
%%%%%%%%%%%%%%%%%%%%%%%%%%%%%%%%%%%%%%%%%%%%%%%%%%%%%%%%%%%%%%%%%%%%%%

%In this section, we present some basic properties on the p.g.f and the moments of the BPVEI, some of which will be useful for the analysis we develop in the following sections.

Having set the notation, in the following result we establish the relations between the p.g.f.s and the moments of the BPVEI and the p.g.f.s and moments of the reproduction and immigration laws, some of which will be useful for the analysis that we develop in the following sections.

\begin{prop}\label{prop:pgf-iter}
The following expressions hold true:
\begin{enumerate}[label=(\alph*),ref=\emph{(\alph*)}]
\item $F_{n+1}(s)=\prod_{i=0}^n h_i(f_{i-1,n}(s))$, for $n\in\N_0$ and $s\in [0,1]$.\label{prop:pgf-iter-a}
\item For $n\in\N_0$,\label{prop:pgf-iter-b}
\begin{enumerate}[label=(\roman*),ref=\emph{(\roman*)}]
\item $E[Z_{n+1}]\hspace*{-0.5ex}=\hspace*{-0.5ex}\sum_{i=0}^n \alpha_{n-i}\prod_{j=0}^i m_{n-j}$.\label{prop:pgf-iter-b-i}
\item $Var[Z_{n+1}]\hspace*{-0.5ex}=\hspace*{-0.5ex}\sum_{i=0}^n\beta_i^2 \prod_{j=i}^n m_j^2+\sum_{i=0}^n\prod_{j=i+1}^n m_j^2\sigma_i^2\left(\alpha_i+\sum_{k=0}^i\alpha_{i-k}\prod_{l=0}^k m_{i-l}\right)$.\label{prop:pgf-iter-b-ii}
\end{enumerate}
\end{enumerate}
\end{prop}

\vspace*{0.5cm}

A BPVEI can be also seen in the following way: each immigrant, when incorporating into the population and independently of the other immigrants, produces a BPVE with initial value the unity. The next proposition formally describes this idea.

\begin{prop}\label{prop:equiv-proc}
Let us consider the independent processes $\{\widetilde{Z}_k^{(j)}\}_{k\in\N_0}$, $j\in\N_0$, defined as:
\begin{equation}\label{eq:def-proc-equiv}
\widetilde{Z}_0^{(j)}=I_j,\quad \widetilde{Z}_{k+1}^{(j)}=\sum_{i=1}^{\widetilde{Z}_k^{(j)}} X_{ki}^{(j)},\quad k\in\N_0,
\end{equation}
where $\{X_{ki}^{(j)}:k\in\N_0;i\in\N;j\in\N_0\}$ is a family of independent random variables such that for each $k\in\N_0$ and $j\in\N_0$ fixed, $X_{ki}^{(j)}$, $i\in\N$ are distributed according to the p.g.f. $f_{k+j}$. Then
\begin{equation*}
Z_n\stackrel{d}{=}\sum_{j=0}^{n-1}  \widetilde{Z}_{n-j}^{(j)},\quad n\in\N,
\end{equation*}
where $\stackrel{d}{=}$ denotes ``equal in distribution''.
\end{prop}

Intuitively, the previous proposition can be interpreted as follows: in each generation $j$, the immigrants $I_j$ generate a BPVE $\{\widetilde{Z}_k^{(j)}\}_{k\in\N_0}$, with environment $v_j=\{f_j,f_{j+1},\ldots\}$ and which starts  with those immigrants $I_j$. Then, the distribution of the number of individuals of the BPVEI at generation $n$ is equal to the distribution of the sum of the total number of descendants in generation $n-j$ of the BPVEs produced by the immigrants in generation $j$, where $j$ goes from the initial generation to the $(n-1)$-st generation. This representation will be useful in Proposition \ref{prop:prob-extinction}, where we establish a necessary and sufficient condition for the certain extinction of the process.

%Thus, the sum , from the initial to the $n$-th generation has the same distribution as  the number of individuals in the generation $n$ of the BPVEI. % and provides an equivalent construction of a BPVEI from a probabilistic point of view.

%%%%%%%%%%%%%%%%%%%%%  EXTINCTION  %%%%%%%%%%%%%%%%%%%%%%%%%%%%%%%%%%
\section{Extinction problem}\label{sec:extinction}
%%%%%%%%%%%%%%%%%%%%%%%%%%%%%%%%%%%%%%%%%%%%%%%%%%%%%%%%%%%%%%%%%%%%%

In this section we address the extinction problem of a BPVEI. It is well known that in the BGWP the incorporation of immigration independent of the generation into the model prevents from the extinction of the population. Similarly, for the BPVEI it is easy to see that 0 is not an absorbing state. Indeed, for each $n\in\N_0$, since $h_n(f_n(0))<1$,
\begin{equation*}
P[Z_{n+1}>0|Z_n=0]=1-P\left[\sum_{i=0}^{I_n}X_{ni}=0\right]=1-\sum_{k=0}^\infty f_n(0)^k P[I_n=k]=1-h_n(f_n(0))>0.
\end{equation*}

However, in this case we can establish some conditions under which the process reaches the state 0 and stays there forever, which can be interpreted as the extinction of the population. To that end, let us denote the probability of this event by
\begin{equation*}\label{def:prob-extinction}
q=P\left[\bigcup_{n=0}^\infty \bigcap_{j=n}^\infty \{Z_j=0\}\right].
\end{equation*}

\begin{prop}\label{prop:prob-extinction}
The probability $q=1$ if and only if
\begin{equation}\label{eq:cond-prob-extinction}
\lim_{n\to\infty}f_{-1,n}(0)=1,\quad\text{ and }\quad \sum_{j=0}^\infty (1-h_j(f_j(0)))<\infty.
\end{equation}
\end{prop}

Note that $f_{-1,n}$ would be the p.g.f. of $Z_{n}$ if $I_0=1$ a.s. and $I_n=0$ a.s. for all $n\in\N$, i.e., if the population starts with the offspring of a single immigrant in the initial generation and there is no immigration after that time. Thus, intuitively, this result indicates that the a.s. extinction of the population happens when the process without immigration becomes extinct with probability one and from some generation, the probability that there are new immigrants in the population and they produce any offspring in the following generations is negligible. This interpretation is illustrated in the examples provided below.

\begin{remark}
\begin{enumerate}[label=(\alph*),ref=\emph{(\alph*)}]
\item It is worthwhile to mention that Proposition \ref{prop:prob-extinction} can be applied even when either the immigration or the offspring distribution are independent of the generation. Indeed, the second condition can also hold true in those cases as we show below, but not for the BGWP with homogeneous immigration.

\item In \cite{kersting-2017}, Theorem 1 provides necessary and sufficient criteria for the first condition in \eqref{eq:cond-prob-extinction} to hold under a mild assumption on the family $\{X_{nj}:\ n\in\N_0;\ j\in\N\}$. In the particular case that the offspring distribution is independent of the generation, it is well known that the aforementioned condition holds if and only if $m\leq 1$, where $m$ is the offspring distribution (see \cite{Athreya-Ney}, p.4).

\item The second requirement in \eqref{eq:cond-prob-extinction} is not difficult to verify. For example, using the mean value theorem, one can prove that a sufficient condition is that the sequence of immigration means is bounded and the series $\sum_{j=0}^\infty (1-f_j(0))$ converges. %Although we do not elaborate them, other necessary or sufficient conditions can be easily obtained by analysing the series $\sum_{j=0}^\infty (1-h_j(f_j(0)))<\infty$.
\end{enumerate}
\end{remark}

\begin{ex}\label{ex:extinct}
\begin{enumerate}[label=(\alph*),ref=\emph{(\alph*)}]
\item Let us consider a BPVEI with $h_n(s)=2^{-1}(1+s)$ for $n=0,1$, and $h_n(s)=1-n^{-2}+n^{-2} s$, $n\geq 2$, whose immigration mean satisfies $\alpha_n=n^{-2}\to 0$.  Observe that for any environment $v=\{f_0,f_1,\ldots\}$, this process fulfils the second condition in Proposition \ref{prop:prob-extinction}. Indeed,\label{ex:extinc-1}
\begin{equation*}
\sum_{n=0}^\infty (1-h_n(f_n(0)))=1-\frac{f_0(0)+f_1(0)}{2}+\sum_{n=2}^\infty \frac{1-f_n(0)}{n^2}\leq 1+\sum_{n=2}^\infty \frac{1}{n^2}<\infty.
\end{equation*}
Thus, if $v$ is such that $f_{-1,n}(0)\to 1$, as $n\to\infty$, then $q=1$ and otherwise, $0\leq q<1$. Note that one can even choose an offspring distribution independent of the generation.

%\begin{figure}[H]
%\centering
%\includegraphics[width=0.3\textwidth]{../sim_ext/graphs/proc_a1a.pdf}\hspace*{0.03\textwidth}
%\includegraphics[width=0.3\textwidth]{../sim_ext/graphs/pnoext_1a.pdf}
%\caption{Left: evolution of the number of individuals (black line) and immigrants (red line) of a BPVEI as described in Example \ref{ex:extinct}~\ref{ex:extinc-1} with $X_{n1}$ following Poisson distribution with mean $m_n=1$, for $n=0,1$, and with mean $m_n=1-n^{-2}$, for $n\geq 2$. Right: evolution of the probability $P[Z_n>0]$ for such a model. }\label{fig:ex-extinction-1a}
%\end{figure}
%
%\begin{figure}[H]
%\centering
%\includegraphics[width=0.3\textwidth]{../sim_ext/graphs/proc_a1b.pdf}\hspace*{0.03\textwidth}
%\includegraphics[width=0.3\textwidth]{../sim_ext/graphs/pnoext_1b.pdf}
%\caption{Left: evolution of the number of individuals (black line) and immigrants (red line) of a BPVEI as described in Example \ref{ex:extinct}~\ref{ex:extinc-1} with $f_n(s)=(7+9s^2)/16$, for $n\in\N_0$. Right: evolution of the probability $P[Z_n>0]$ for such a model. }\label{fig:ex-extinction-1b}
%\end{figure}

\item Let us consider a BPVEI with $h_n(s)=2^{-1}(1+s)$, whose immigration mean is $\alpha_n=2^{-1}$, $n\in\N_0$ and $f_n(s)=2^{-1}(s+1)$, for $n=0,1$, and $f_n(s)=1-n^{-2}+n^{-2} s$, for $n\geq 2$. Note that in this case, the immigration is independent of the generation and the iteration of the offspring p.g.f.s leads to $f_{-1,n}(s)=1-4^{-1}(n!)^{-2}+4^{-1}(n!)^{-2}s$, $n\geq 1$, from which one obtains $f_{-1,n}(0)\to 1$, as $n\to\infty$. Moreover, the process also satisfies the second condition in Proposition \ref{prop:prob-extinction}, as we show below, and as a consequence, $q=1$. \label{ex:extinc-2}
\begin{equation*}
\sum_{n=0}^\infty (1-h_n(f_n(0)))=\sum_{n=0}^\infty \frac{1-f_n(0)}{2}=\frac{1}{2}+\sum_{n=2}^\infty \frac{1}{2n^2}<\infty.
\end{equation*}

%\begin{figure}[H]
%\centering
%\includegraphics[width=0.3\textwidth]{../sim_ext/graphs/proc_a2.pdf}\hspace*{0.03\textwidth}
%\includegraphics[width=0.3\textwidth]{../sim_ext/graphs/pnoext_2.pdf}\hspace*{0.03\textwidth}
%\caption{Left: evolution of the number of individuals (black line) and immigrants (red line) of a BPVEI as described in Example \ref{ex:extinct}~\ref{ex:extinc-2}. Right: evolution of the probability $P[Z_n>0]$ for a such model.}\label{fig:ex-extinction-2}
%\end{figure}

\item Let us consider a BPVEI with the same immigration laws as in the previous example but offspring distributions $f_n(s)=2^{-1}(s+1)$, for $n=0,1$, and $f_n(s)=1-n^{-1}+n^{-1} s$, for $n\geq 2$. With the same arguments as before, one can show that $0\leq q<1$ since $f_{-1,n}(s)=1-4^{-1}(n!)^{-1}+4^{-1}(n!)^{-1}s$, $n\geq 1$, implying $f_{-1,n}(0)\to 1$, as $n\to\infty$, and the second condition in Proposition \ref{prop:prob-extinction} does not hold true due to the fact that \label{ex:extinc-3}
\begin{equation*}
\sum_{n=0}^\infty (1-h_n(f_n(0)))=\sum_{n=0}^\infty \frac{1-f_n(0)}{2}=\frac{1}{2}+\sum_{n=2}^\infty \frac{1}{2n}=\infty.
\end{equation*}

%\begin{figure}[H]
%\centering
%\includegraphics[width=0.3\textwidth]{../sim_ext/graphs/proc_a3.pdf}\hspace*{0.03\textwidth}
%\includegraphics[width=0.3\textwidth]{../sim_ext/graphs/pnoext_3.pdf}\hspace*{0.03\textwidth}
%\caption{Left: evolution of the number of individuals (black line) and immigrants (red line) of a BPVEI as described in Example \ref{ex:extinct}~\ref{ex:extinc-3}. Right: evolution of the probability $P[Z_n>0]$ for a such model.}\label{fig:ex-extinction-3}
%\end{figure}

\begin{figure}[H]
\centering
\includegraphics[width=0.45\textwidth]{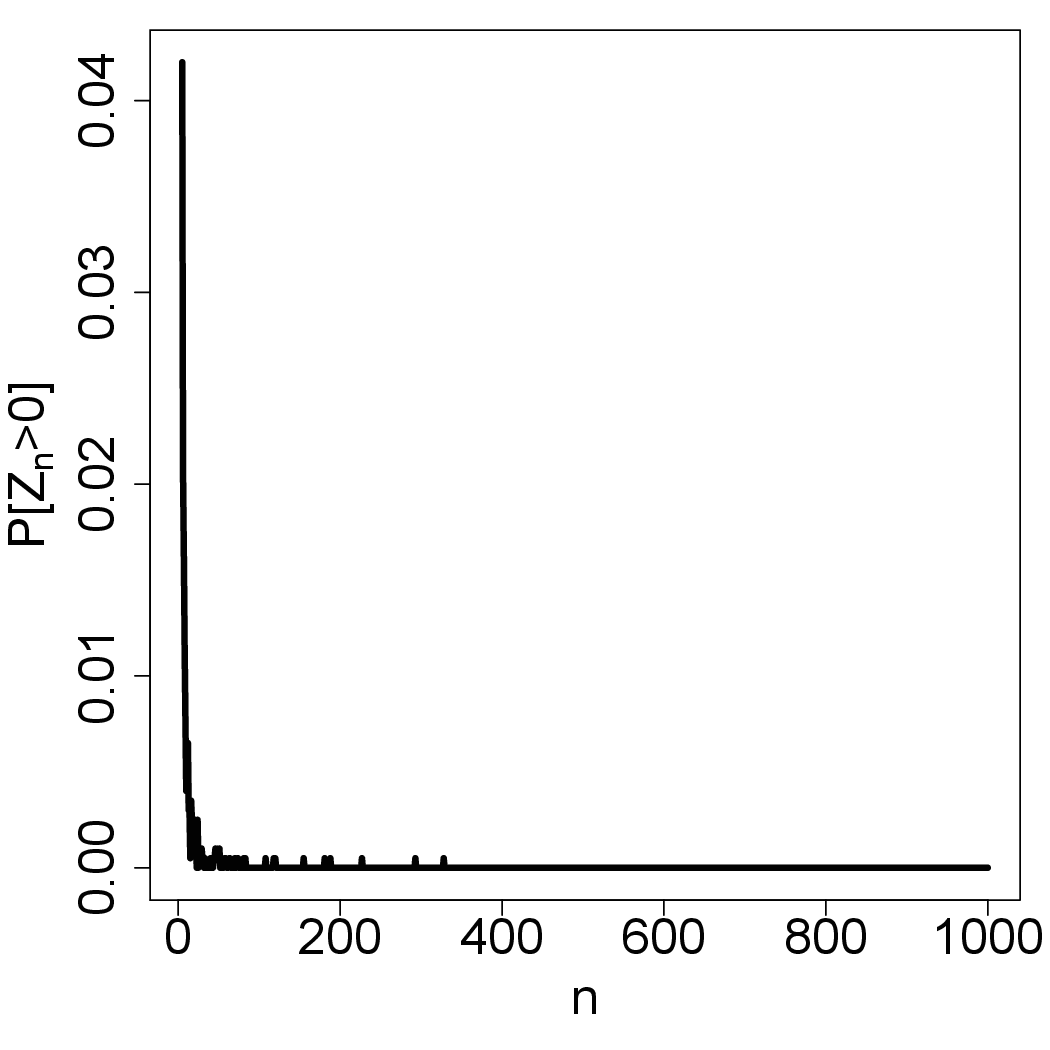}\hspace*{0.03\textwidth}
\includegraphics[width=0.45\textwidth]{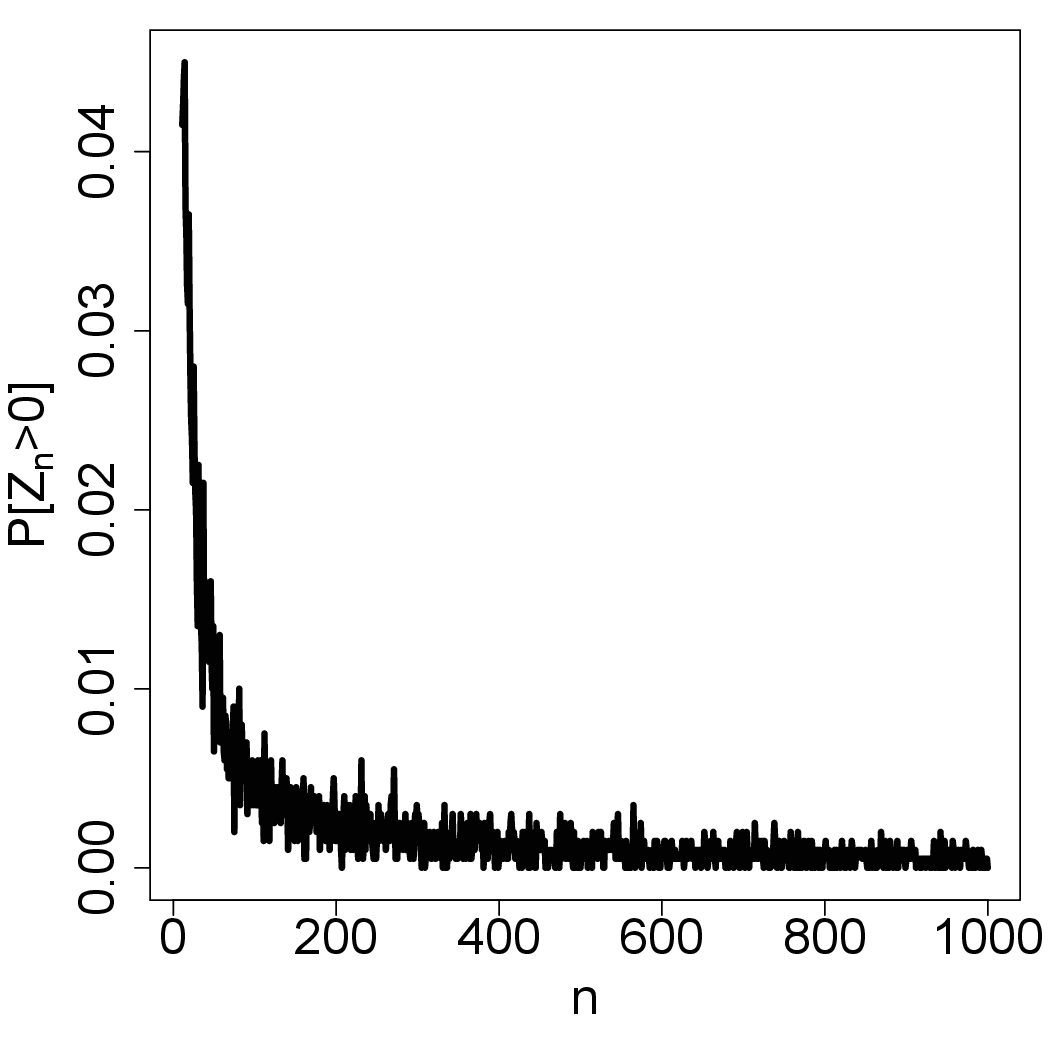}\hspace*{0.03\textwidth}
\caption{Comparison of the evolution of the probability $P[Z_n>0]$ for the models described in Examples \ref{ex:extinct}~\ref{ex:extinc-2} (left) and \ref{ex:extinct}~\ref{ex:extinc-3} (right). The probabilities $P[Z_n>0]$, $n\in\N$, were estimated by Monte Carlo methods at each generation $n\in\N$. To that end, the first 1000 generations of 2000 processes following the models in the aforementioned examples were simulated using the statistical software and programming environment \texttt{R} (see \cite{R}).}\label{fig:ex-extinction}
\end{figure}
\end{enumerate}
\end{ex}

%%%%%%%%%%%%%%%%%%% LIMIT THEOREMS %%%%%%%%%%%%%%%%%%%%%%%%%%%%%%
\section{Asymptotic behaviour}\label{sec:lim-behaviour}
%%%%%%%%%%%%%%%%%%%%%%%%%%%%%%%%%%%%%%%%%%%%%%%%%%%%%%%%%%%%%%%%%%%%%

In this section, we analyse the limiting behaviour of BPVEIs focusing our attention on those with critical offspring distributions according to the classification established in \cite{kersting-2017} for BPVEs, that is, with reproduction laws satisfying
\begin{equation}\label{eq:critical-offs}
\sum_{k=0}^n \frac{\nu_k}{\mu_{k-1}}\to\infty\quad\text{ and }\quad \frac{1}{\mu_n}=o\left(\sum_{k=0}^n \frac{\nu_k}{\mu_{k-1}}\right),\quad \text{ as }n\to\infty.
\end{equation}
Note that in \eqref{eq:critical-offs} the convergence of the sequence $\{\mu_n\}_{n\in\N_0}$ is not required; the limit may not exist, or it may be 0, infinity, or any positive value.

\vspace*{2ex}

Moreover, it is important to remark that in contrast to the BGWP, for BPVEs one might encounter different growth rates (see \cite{MaCphee-Schuh-1983}). In order to avoid that situation, we consider the following regularity condition introduced in \cite{kersting-2017} for the BPVE and which is satisfied by a wide class of probability distributions: for every $\epsilon>0$ there is a constant $c_\epsilon<\infty$ such that for all $n\in\N_0$,
\begin{equation}\label{eq:cond-regularity}
E[X_{n1}^2;X_{n1}>c_\epsilon(1+E[X_{n1}])]\leq \epsilon E[X_{n1}^2;X_{n1}\geq 2].
\end{equation}
%This uniformity condition   %as explained in \cite{kersting-2017}, p.9.

\vspace*{2ex}
%In this section, we analyse the limiting behaviour of BPVEIs focusing our attention on those with offspring distributions which are critical according to the classification established in \cite{kersting-2017} for BPVEs, that is, with reproduction laws satisfying
%\begin{equation}\label{eq:critical-offs}
%\sum_{k=0}^n \frac{\nu_k}{\mu_{k-1}}\to\infty\quad\text{ and }\quad \frac{1}{\mu_n}=o\left(\sum_{k=0}^n \frac{\nu_k}{\mu_{k-1}}\right),\quad \text{ as }n\to\infty.
%\end{equation}
%Note that in \eqref{eq:critical-offs} the convergence of the sequence $\{\mu_n\}_{n\in\N_0}$ is not required; the limit may not exist, or it may be 0, infinity, or any positive value.
%
%\vspace*{2ex}
%
%
%Moreover, we consider the following regularity condition: for every $\epsilon>0$ there is a constant $c_\epsilon<\infty$ such that for all $n\in\N_0$,
%\begin{equation}\label{eq:cond-regularity}
%E[X_{n1}^2;X_{n1}>c_\epsilon(1+E[X_{n1}])]\leq \epsilon E[X_{n1}^2;X_{n1}\geq 2].
%\end{equation}
%This condition is satisfied by a wide class of probability distributions. It was introduced in \cite{kersting-2017} for the BPVE and prevents it from having different growth rates. %as explained in \cite{kersting-2017}, p.9.
%
%\vspace*{2ex}

The following result is a counterpart of Theorem VI.7.4 in \cite{Athreya-Ney} for BGWPs with immigration independent of the generation and establishes the asymptotic distribution of a BPVEI suitably normalized when the immigration means stabilize.

\begin{thm}\label{thm:asym-critical}
Let $\{Z_n\}_{n\in\N_0}$ be a BPVEI satisfying \eqref{eq:critical-offs} and \eqref{eq:cond-regularity} and denote
\begin{equation*}
a_{n+1}=\frac{\mu_n}{2}\sum_{k=0}^n \frac{\nu_k}{\mu_{k-1}},\quad n\in\N_0.
\end{equation*}
%Assume that $\nu_n\to\nu>0$ and $\alpha_n\to\alpha>0$, as $n\to\infty$, $0<\tau=\inf_{n\in\N_0} h_n(0)$ and $\sup_{n\in\N_0} h_n''(1)<\infty$.
Assume that $\nu_n\to\nu>0$ and $\alpha_n\to\alpha>0$, as $n\to\infty$, $0<\tau=\inf_{n\in\N_0} h_n(0)$ and $\sup_{n\in\N_0} \beta_n^2<\infty$.
Then, the asymptotic distribution of $Z_n/a_{n}$ is a gamma distribution with parameters $2\alpha/\nu$ and 1.
\end{thm}

\begin{cor}\label{cor:prob-non-ext}
Let $\{Z_n\}_{n\in\N_0}$ be a BPVEI satisfying the conditions of Theorem \ref{thm:asym-critical}, then
\begin{equation*}
P[Z_n>0]\to 1,\quad\text{ as }n\to\infty.
\end{equation*}
\end{cor}

\vspace*{2ex}

Now, we shall provide a description of the relation between the asymptotic behaviour of the BPVEI given by Theorem \ref{thm:asym-critical} and the asymptotic behaviour of the BPVEs $\{\widetilde{Z}_k^{(j)}\}_{n\in\N_0}$ introduced in   Proposition \ref{prop:equiv-proc}. First, note that if $\{Z_n\}_{n\in\N_0}$ is a BPVEI with critical offspring distribution, then the processes $\{\widetilde{Z}_k^{(j)}\}_{n\in\N_0}$ are critical BPVEs. On the other hand, at each generation $j$, $\alpha_j$ individuals are expected to immigrate to the population. Since from Theorem 4 in \cite{kersting-2017} one has that the probability that the BPVE produced by one immigrant at generation $j$, $\{\bar{Z}_n^{(j)}\}_{n\in\N_0}$, is not extinct after $n$ generations is $P[\overline{Z}_n^{(j)}>0]\sim 2\left(\sum_{k=j}^{j+n} \frac{\nu_k}{\mu_{k-1}}\right)^{-1}$, as $n\to\infty$, then one obtains that $P[\widetilde{Z}_n^{(j)}>0]\sim 2\alpha_j\left(\sum_{k=j}^{j+n} \frac{\nu_k}{\mu_{k-1}}\right)^{-1}$. Moreover, Theorem 4 in \cite{kersting-2017} provides the asymptotic distribution of those BPVEs; indeed, for each $j\in\N_0$, each one of these BPVEs satisfies that the distribution of the process $\bar{Z}_n^{(j)}/a_n^{(j)}$ conditioned on the event $\{\bar{Z}_n^{(j)}>0\}$ converges to an exponential distribution with mean equal to one, where 
$$a_n^{(j)}=\frac{\mu_{j+n}}{2\mu_{j-1}}\sum_{k=j}^{j+n} \frac{\nu_k}{\mu_{k-1}},\quad n\geq j.$$
%Now, we shall provide a description of the relation between the asymptotic behaviour of the BPVEI given by Theorem \ref{thm:asym-critical} and the asymptotic behaviour of the BPVEs $\{\widetilde{Z}_k^{(j)}\}_{n\in\N_0}$ introduced in   Proposition \ref{prop:equiv-proc}. First, note that if $\{Z_n\}_{n\in\N_0}$ is a BPVEI with critical offspring distribution, then the processes $\{\widetilde{Z}_k^{(j)}\}_{n\in\N_0}$ are critical BPVEs. On the other hand, at each generation $j$, $\alpha_j$ individuals are expected to immigrate to the population. Since from Theorem 4 in \cite{kersting-2017} one has that the probability that the BPVE produced by one immigrant at generation $j$, $\{\bar{Z}_n^{(j)}\}_{n\in\N_0}$, is not extinct after $n$ generations is $P[\overline{Z}_n^{(j)}>0]\sim 2\left(\sum_{k=j}^{j+n} \frac{\nu_k}{\mu_{k-1}}\right)^{-1}$, as $n\to\infty$, then one obtains that the BPVE generated by all the immigrants at the $j$-th generation is non-extinct after $n$ generations with a probability of the order $2\alpha_j\left(\sum_{k=j}^{j+n} \frac{\nu_k}{\mu_{k-1}}\right)^{-1}$. Moreover, Theorem 4 in \cite{kersting-2017} provides the asymptotic distribution of those BPVEs; indeed, for each $j\in\N_0$, each one of these BPVEs satisfies that $\bar{Z}_n^{(j)}/a_n^{(j)}$ converges to an exponential distribution with mean equal to one, where $$a_n^{(j)}=\frac{\mu_{n+j}}{2\mu_{j-1}}\sum_{k=j}^{j+n} \frac{\nu_k}{\mu_{k-1}},\quad n\geq j.$$

\vspace*{2ex}

We also remark here that Theorem \ref{thm:asym-critical} can be applied even in the cases when one of the distributions of the model -either the offspring distribution or the immigration law- is inhomogeneous but the other one is homogeneous.
%the offspring distribution is homogeneous but the immigration is not homogeneous and when the offspring distribution is inhomogeneous and the immigration law is homogeneous. ç
However, Theorem \ref{thm:asym-critical} cannot be applied for the model considered in \cite{Mitov-Omey-2014} due to the fact that in that situation $\nu_n\to\infty$, as $n\to\infty$.

\section*{Appendix}
%%%%%%%%%%%%%%%%%%%%%%%%%%%%%%%%%%%%%%%%%%%%%%%%%%%%%%%%%%%%%%%%%%%%%

\begin{Prf}[Proposition \ref{prop:pgf-iter}]
\ref{prop:pgf-iter-a} Using the independence of the offspring variables and the independence between the offspring and immigration phases, one has
\begin{equation}\label{eq:pgf}
F_{n+1}(s)=E\left[E\left[s^{\sum_{i=1}^{Z_n+I_n}X_{ni}}|Z_n+I_n\right]\right]=E[f_n(s)^{Z_n+I_n}]=F_n(f_n(s))h_n(f_n(s)),
\end{equation}
and by iterating the previous formula one obtains the  expression in \ref{prop:pgf-iter-a}.

\ref{prop:pgf-iter-b}~\ref{prop:pgf-iter-b-i} It can be easily obtained using a recursion on $E[Z_{n+1}]=m_n(E[Z_n]+\alpha_n)$, which follows from
$$E[Z_{n+1}|Z_n]=E\left[\sum_{i=1}^{Z_n+I_n}X_{ni}|Z_n\right]=m_n(\alpha_n+Z_n)\quad a.s.$$

\ref{prop:pgf-iter-b}~\ref{prop:pgf-iter-b-ii} To that end, we make use of the fact that $Var[Z_{n+1}]=F_{n+1}''(1)+F_{n+1}'(1)-F_{n+1}'(1)^2$ and it is straightforward from \eqref{eq:pgf} that
\begin{align*}
F_{n+1}'(1)&=F_n'(1)f_n'(1)+h_n'(1)f_n'(1),\\
F_{n+1}''(1)&=F_n''(1)f_n'(1)^2+F_n'(1)f_n''(1)+2 F_n'(1)f_n'(1)^2h_n'(1)+f_n'(1)^2h_n''(1)+f_n''(1)h_n'(1).
\end{align*}
An arrangement of the terms leads to
\begin{align*}
Var[Z_{n+1}]&=f_n'(1)^2[F_n''(1)+F_n'(1)-F_n'(1)^2]+F_n'(1)[f_n''(1)+f_n'(1)-f_n'(1)^2]\\
&\phantom{=}+f_n'(1)^2[h_n''(1)+h_n'(1)-h_n'(1)^2]+h_n'(1)[f_n''(1)+f_n'(1)-f_n'(1)^2]\\
&=m_n^2 (Var [Z_n]+\beta_n^2)+\sigma_n^2 (E[Z_n]+\alpha_n),
\end{align*}
and to obtain the desired formula we use a recursion in the  following expression
\begin{equation*}
\frac{Var[Z_{n+1}]}{\prod_{j=0}^n m_j^2}=\frac{Var[Z_n]}{\prod_{j=0}^{n-1}m_j^2}+\frac{\beta_n^2m_n^2+\sigma_n^2 (E[Z_n]+\alpha_n)}{\prod_{j=0}^n m_j^2}.%=\sum_{i=0}^n \frac{\beta_i^2m_i^2+\sigma_i^2 (E[Z_i]+\alpha_i)}{\prod_{j=0}^i m_j^2}.
\end{equation*}
\end{Prf}

%\begin{Prf}[Proposition \ref{prop:pgf}]
%\ref{prop:pgf-a} It is straightforward by iterating the formula in Proposition \ref{prop:pgf-iter}~\ref{prop:pgf-iter-a}.
%\end{Prf}

\vspace*{0.75cm}

\begin{Prf}[Proposition \ref{prop:equiv-proc}]
Let $\widetilde{F}_n^{(j)}$ be the p.g.f. of $\widetilde{Z}_n^{(j)}$, $n,j\in\N_0$. Then, for $s\in [0,1]$ and $j\in\N_0$, $\widetilde{F}_0^{(j)}(s)=h_j(s)$ and for each $k\in\N_0$,
\begin{equation*}
\widetilde{F}_{k+1}^{(j)}(s)=E\left[E\left[s^{\widetilde{Z}_{k+1}^{(j)}}|\widetilde{Z}_k^{(j)}\right]\right]=E\left[\left(f_{k+j}(s)\right)^{\widetilde{Z}_k^{(j)}}\right]=\widetilde{F}_{k}^{(j)}(f_{k+j}(s))=h_j(f_{j-1,j+k}(s)).
\end{equation*}
As a consequence, since the processes $\{\widetilde{Z}_{k}^{(j)}\}_{k\in\N_0}$, $j\in\N_0$, are independent, one has
\begin{equation*}
E\left[s^{\sum_{j=0}^{n-1}\widetilde{Z}_{n-j}^{(j)}}\right]=\prod_{j=0}^{n-1} E\left[s^{\widetilde{Z}_{n-j}^{(j)}}\right]=\prod_{j=0}^{n-1} h_j(f_{j-1,n-1}(s)) =F_n(s).
\end{equation*}
\end{Prf}

\begin{Prf}[Proposition \ref{prop:prob-extinction}]
First, one has $q=\lim_{n\to\infty} P\left[\cap_{j=n}^\infty\{Z_j=0\}\right]$ and
\begin{equation*}
P\left[\cap_{j=n}^\infty\{Z_j=0\}\right]=P[Z_n=0]\prod_{j=n}^\infty P\left[Z_{j+1}=0|Z_j=0\right]=  F_n(0)\prod_{j=n}^\infty h_j(f_j(0)).
\end{equation*}
Thus, a necessary and sufficient condition for $q=1$ is that both limits $F_n(0)\to 1$ and $\prod_{j=n}^\infty h_j(f_j(0))\to 1$, as $n\to\infty$. Now, we prove that this condition is equivalent to that established in Proposition \ref{prop:prob-extinction}.

Note that $\prod_{j=n}^\infty h_j(f_j(0))\to 1$, as $n\to\infty$ holds if and only if $\sum_{j=n}^\infty \log (h_j(f_j(0)))\to 0$, as $n\to\infty$, which is equivalent to the condition $\sum_{j=0}^\infty (1-h_j(f_j(0)))<\infty$.

Regarding the first limit, from Proposition \ref{prop:pgf-iter}~\ref{prop:pgf-iter-a} one has that if $F_n(0)\to 1$, then \linebreak$h_i(f_{i-1,n}(0))\to 1$, for each $i=0,1,\ldots,n$, which implies $f_{-1,n}(0)\to 1$, as $n\to\infty$.

Conversely, assume that $f_{-1,n}(0)\to 1$, as $n\to\infty$, and $\sum_{j=0}^\infty (1-h_j(f_j(0)))<\infty$. Observe that if we prove $P\left[Z_n>0\ \text{i.o.}\right]=0$, writing ``i.o.'' to mean infinitely often, then $F_n(0)\to 1$, as $n\to\infty$.  To that end, we use the representation of the BPVEI described in Proposition \ref{prop:equiv-proc} and define the random variable $L=\max\{j\in\N:\widetilde{Z}_1^{(j)}>0\}$. Since $P[\widetilde{Z}_1^{(j)}>0]=1-h_j(f_j(0))$, $j\in\N_0$, from Borel-Cantelli's lemma,
$$P\left[\widetilde{Z}_1^{(j)}>0\text{ f.o.}\right]=1,$$
where ``f.o.'' stands for finitely often, and $L<\infty$ a.s. Consequently,
\begin{align*}
P\left[Z_n>0\ \text{i.o.}\right]&=P\left[Z_n>0\ \text{i.o.},L<\infty\right]\\
&=\sum_{k=0}^\infty P\left[Z_n>0\ \text{i.o.},L=k\right]\\
&\leq\sum_{k=0}^\infty \sum_{j=0}^k P\left[\widetilde{Z}_n^{(j)}>0\ \text{i.o.}\right].
\end{align*}
and for each $j\in\N_0$ fixed,
\begin{align*}
P\left[\widetilde{Z}_n^{(j)}>0\ \text{i.o.}\right]&=1-P\left[\widetilde{Z}_n^{(j)}=0\ \text{f.o.}\right]\\
&=1-\lim_{n\to\infty} \widetilde{F}_n^{(j)}(0)\\
&=1-\lim_{n\to\infty} h_j(f_{j-1,j+n-1}(0)),
\end{align*}
where recall that $\widetilde{F}_{n}^{(j)}(s)=h_j(f_{j-1,j+n-1}(s))$ is the p.g.f. of $\widetilde{Z}_n^{(j)}$, $n\in\N$. Now, since $f_{-1,n}(0)\to 1$, as $n\to\infty$, $f_{j,n}(0)\to 1$ and $h_j(f_{j-1,n}(0))\to 1$, for each $j\in\N_0$, which implies $P\left[Z_n>0\ \text{i.o.}\right]=0$ and hence, $F_n(0)\to 1$.
\end{Prf}

\vspace*{0.75cm}

In order to prove Theorem \ref{thm:asym-critical}, we make use of the shape functions corresponding to the p.g.f.s of the reproduction laws. The \emph{shape function} of the p.g.f. $f_k$, $k\in\N_0$, is the function $\varphi_k:[0,1)\rightarrow \R$ satisfying
\begin{equation}\label{eq:shape-function}
\frac{1}{1-f_k(s)}=\frac{1}{m_k(1-s)}+\varphi_k(s),\quad s\in [0,1).
\end{equation}
This function satisfies $\varphi_k(s)\geq 0$, for $0\leq s<1$, and can be extended by setting $
\varphi_k(1)=\nu_k/2$. The reader is referred to Section 3 in \cite{kersting-2017} for further details about the properties of these functions. Then, using the fact that $f_{k,n}=f_{k+1}\circ f_{k+1,n}$, for $k=0,\ldots,n-1$, by iterating \eqref{eq:shape-function} one obtains:
\begin{equation}\label{eq:iter-shap}
\frac{1}{1-f_{k,n}(s)}=\frac{\mu_k}{\mu_n(1-s)}+\mu_k\sum_{l=k+1}^n \frac{\varphi_l(f_{l,n}(s))}{\mu_{l-1}},\quad s\in [0,1).
\end{equation}
To prove Theorem \ref{thm:asym-critical}, we also need the next lemma; its proof is similar to the one of Lemma 7 in \cite{kersting-2017} and is omitted.

\begin{lem}\label{lem:conv-sum-phi}
Let $i=0,\ldots,n$ be fixed. Let $\{Z_n\}_{n\in\N_0}$ be a BPVEI satisfying \eqref{eq:critical-offs} and \eqref{eq:cond-regularity}, then
\begin{equation*}
\sup_{s\in [0,1]} \left|\sum_{k=i}^n \frac{\varphi_k(f_{k,n}(s))}{\mu_{k-1}}-\sum_{k=i}^n \frac{\varphi_k(1)}{\mu_{k-1}}\right|=o\left(\sum_{k=i}^n \frac{\varphi_k(1)}{\mu_{k-1}}\right),\quad \text{ as }n\to\infty.
\end{equation*}
\end{lem}

\vspace*{0.75cm}

\begin{Prf}[Theorem \ref{thm:asym-critical}]
We shall prove that the Laplace transform of $Z_n/a_n$ converges to the Laplace transform of a gamma distribution with parameters $2\alpha/\nu$ and 1. Let us fix $\lambda> 0$ and for simplicity, let us denote $s_n=e^{-\lambda/a_n}$, $n\in\N$. From Proposition \ref{prop:pgf-iter}~\ref{prop:pgf-iter-a} and using a Taylor expansion for the functions $\log h_i$ at 1, for $i=0,\ldots,n$, one has
\begin{align*}
F_{n+1}(s_{n+1})&=\exp\bigg\{\sum_{i=0}^n \log(h_i(f_{i-1,n}(s_{n+1})))\bigg\}\nonumber\\
&=\exp\bigg\{-\sum_{i=0}^n \alpha_i(1-f_{i-1,n}(s_{n+1}))\\
&\phantom{=}+\frac{1}{2}\sum_{i=0}^n \left(\frac{h_i''(\xi_{in})}{h_i(\xi_{in})}-\frac{h_i'(\xi_{in})^2}{h_i(\xi_{in})^2}\right)(1-f_{i-1,n}(s_{n+1}))^2\bigg\},
\end{align*}
with $f_{i-1,n}(s_{n+1})< \xi_{in}< 1$, $i=0,\ldots,n$, $n\in\N_0$. Thus, the result yields by proving the following convergences, as $n\to\infty$,
\begin{align}
\sum_{i=0}^n \alpha_i(1-f_{i-1,n}(s_{n+1}))&\to \log(1+\lambda)^{\frac{2\alpha}{\nu}},\label{eq:result1-thm-asym-critical}\\
\sum_{i=0}^n \left(\frac{h_i''(\xi_{in})}{h_i(\xi_{in})}-\frac{h_i'(\xi_{in})^2}{h_i(\xi_{in})^2}\right)(1-f_{i-1,n}(s_{n+1}))^2&\to 0.\label{eq:result2-thm-asym-critical}
\end{align}

We shall start with \eqref{eq:result1-thm-asym-critical}. To prove this result, for each $N,n\in\N_0$, $n\geq N$, let  us denote
\begin{equation*}
S_{1n}^{(N)}=\sum_{i=0}^{N-1}\alpha_i(1-f_{i-1,n}(s_{n+1})),\quad \text{ and }\quad S_{2n}^{(N)}=\sum_{i=N}^{n}\alpha_i(1-f_{i-1,n}(s_{n+1})).
\end{equation*}
By applying \eqref{eq:iter-shap},
\begin{align*}
S_{2n}^{(N)}&=\sum_{i=N}^{n}\frac{\alpha_i}{\frac{\mu_{i-1}}{\mu_n(1-s_{n+1})}+\mu_{i-1}\sum_{l=i}^n \frac{\varphi_l(f_{l,n}(s_{n+1}))}{\mu_{l-1}}}\\
&=\sum_{i=N}^{n}\frac{\alpha_i}{\varphi_i(f_{i,n}(s_{n+1}))}\cdot\frac{\varphi_i(f_{i,n}(s_{n+1}))}{\mu_{i-1}}\cdot\frac{1}{\frac{1}{\mu_n(1-s_{n+1})}+\sum_{l=i}^n \frac{\varphi_l(f_{l,n}(s_{n+1}))}{\mu_{l-1}}},
\end{align*}
where the functions $\varphi_l(\cdot)$, $l=0,\ldots,n$, were introduced in \eqref{eq:shape-function}. First, we deal with this latter term.

Let $0<\eps<\frac{2\alpha}{\nu}$. We shall prove that there exist $J_0=J_{0,\eps},N_0=N_{0,\eps}\in\N$ such that $N_0>J_0$ and for  $n\geq N_0$ and $J_0\leq i\leq n$,
\begin{equation}\label{eq:proof-asymp-ineq}
\left|\frac{\alpha_i}{\varphi_i(f_{i,n}(s_{n+1}))}-\frac{2\alpha}{\nu}\right|<\eps.
\end{equation}

From the proof of Lemma 7 in \cite{kersting-2017}, one has that given $\varrho>0$, there exists $0<\eta=\eta_\varrho<1$, such that
\begin{equation}\label{eq:eta-ti-lemma}
\sup_{t_i\leq t\leq 1}|\varphi_i(1)-\varphi_i(t)|<\varrho \varphi_i(1),
\end{equation}
for $t_i=1-\frac{\eta}{1+m_i}$, $i\in\N_0$. To obtain \eqref{eq:proof-asymp-ineq}, we prove that given $\varrho>0$, there exists %$N_1=N_{1,\varrho},I_1=I_{1,\varrho}\in\N$,  $N_1>I_1$
$N_1=N_{1,\varrho}$, such that $f_{i,n}(s_{n+1})\geq t_i$, %for $I_1\leq i\leq n$, $n\geq N_1$, or equivalently,
for $0\leq i\leq n$, $n\geq N_1$, or equivalently,
%\begin{equation}\label{eq:proof-inf-sup}
%\inf_{n\geq N_1}\inf_{I_1\leq i\leq n} f_{i,n}(s_{n+1})\geq \sup_{n\geq N_1} \sup_{I_1\leq i\leq n} t_i.
%\end{equation}
%\begin{equation}\label{eq:proof-inf-sup}
%f_{i,n}(s_{n+1})\geq t_i,\quad \text{ for }I_1\leq i\leq n, n\geq N_1.
%\end{equation}
\begin{equation}\label{eq:proof-inf-sup}
(1-f_{i,n}(s_{n+1}))(1+m_i)\leq\eta,\quad \text{ for }0\leq i\leq n, n\geq N_1.
\end{equation}

We write $s_{n+1}=1-\frac{\lambda}{a_{n+1}}+\frac{y_n}{a_{n+1}}$, with $y_n\to 0$; using $f_{i,n}'(1)=\prod_{j=i+1}^n f_j'(1)=\frac{\mu_n}{\mu_i}$ and a Taylor expansion for $f_{i,n}(\cdot)$ at 1, we have
\begin{align*}
(1-f_{i,n}(s_{n+1}))(1+m_i)&\leq f_{i,n}'(1)\left(1-s_{n+1}\right)(1+m_i)\\
&\leq \frac{\mu_n}{\mu_i}\left(\frac{\lambda}{a_{n+1}}-\frac{y_n}{a_{n+1}}\right)(1+m_i)\\
&\leq \frac{\mu_n}{a_{n+1}}(\lambda-y_n)\left(\frac{1}{\mu_i}+\frac{1}{\mu_{i-1}}\right)\\
&\leq (\lambda-y_n)\left(\frac{1}{a_{i+1}}+\frac{1}{a_{i}}\right).
\end{align*}
Now, since $a_i\to\infty$, as $i\to\infty$, given $\epsilon>0$, there exists $J_1=J_{1,\epsilon}$ such that $a_i^{-1}\leq \epsilon(2\lambda)^{-1}$, for $i\geq J_1$, and consequently
\begin{equation*}
(1-f_{i,n}(s_{n+1}))(1+m_i)\leq\frac{\epsilon}{\lambda}(\lambda-y_n) ,\quad n\geq J_1,\ J_1\leq i\leq n,
\end{equation*}
and therefore,
\begin{equation}\label{eq:bound-proof-lim}
\lim_{n\to\infty}\sup_{J_1\leq i\leq n}(1-f_{i,n}(s_{n+1}))(1+m_i)\leq \epsilon.
\end{equation}
On the other hand,
\begin{equation*}
\sup_{0\leq i\leq J_1-1}(1-f_{i,n}(s_{n+1}))(1+m_i)\leq \frac{\mu_n}{a_{n+1}}(\lambda-y_n)\max_{0\leq i\leq J_1-1}\left(\frac{1}{\mu_i}+\frac{1}{\mu_{i-1}}\right),
\end{equation*}
and hence,
\begin{equation*}
\lim_{n\to\infty} \sup_{0\leq i\leq J_1-1}(1-f_{i,n}(s_{n+1}))(1+m_i)=0,
\end{equation*}
which with \eqref{eq:bound-proof-lim} and taking into account that $\epsilon$ is arbitrary gives us
\begin{equation}\label{eq:lim-sup-proof-asymp}
\lim_{n\to\infty} \sup_{0\leq i\leq n}(1-f_{i,n}(s_{n+1}))(1+m_i)=0.
\end{equation}
From this, \eqref{eq:proof-inf-sup} yields. Moreover, since $\varphi_i(1)\to\nu/2>0$, as $i\to\infty$, there exists $J_2\in\N$ such that $\varphi_i(1)>0$, for each $i\geq J_2$, and we can safely assume that $N_1\geq J_2$. Taking these facts into account together with \eqref{eq:eta-ti-lemma} and \eqref{eq:proof-inf-sup}, one obtains
\begin{equation}\label{eq:unif-conv-phi}
\sup_{n\geq N_1}\sup_{J_2\leq i\leq n} \left|\frac{\varphi_i(f_{i,n}(s_{n+1}))}{\varphi_i(1)}-1\right|\leq \sup_{n\geq N_1} \sup_{J_2\leq i\leq n} \sup_{t_i\leq t\leq 1}\left|\frac{\varphi_i(t)}{\varphi_i(1)}-1\right|<\varrho.
\end{equation}

Now, the claim in \eqref{eq:proof-asymp-ineq} yields from \eqref{eq:unif-conv-phi}, the convergences $\alpha_i\to\alpha$ and $\varphi_i(1)=\nu_i/2\to \nu/2$, as $i\to\infty$, and the inequality
\begin{equation*}
\left| \frac{\alpha_i}{\varphi_i(f_{i,n}(s_{n+1}))}-\frac{2\alpha}{\nu}\right|\leq \frac{\alpha_i}{\varphi_i(1)}\left| \frac{\varphi_i(1)}{\varphi_i(f_{i,n}(s_{n+1}))}-1\right|+\left|\frac{2\alpha_i}{\nu_i}-\frac{2\alpha}{\nu}\right|.
\end{equation*}

Let us denote
\begin{equation*}\label{eq:partition-sn}
x_{j,n}=\frac{1}{\mu_n(1-s_{n+1})}+\sum_{l=n-j}^n \frac{\varphi_l(f_{l,n}(s_{n+1}))}{\mu_{l-1}},\quad j=-1,0,\ldots,n;n\in\N_0,
\end{equation*}
where the empty sum is considered to be 0; then, $b_n=x_{-1,n},x_{0,n},\ldots,x_{n-1,n},x_{n,n}=d_n$ is a partition of the interval $[b_n,d_n]$ and
\begin{equation*}
\sum_{i=N}^{n}\frac{\varphi_i(f_{i,n}(s_{n+1}))}{\mu_{i-1}}\cdot\frac{1}{\frac{1}{\mu_n(1-s_{n+1})}+\sum_{l=i}^n \frac{\varphi_l(f_{l,n}(s_{n+1}))}{\mu_{l-1}}}=\sum_{i=N}^{n} x_{n-i,n}^{-1}(x_{n-i,n}-x_{n-i-1,n}).
\end{equation*}

Combining this with \eqref{eq:proof-asymp-ineq}, for $N\geq J_0$, and $n\geq \max\{N_0,N\}$ one has
\begin{equation}\label{eq:bounds-S2n}
\left(\frac{2\alpha}{\nu}-\eps\right)\sum_{i=N}^{n} x_{n-i,n}^{-1}(x_{n-i,n}-x_{n-i-1,n})\leq S_{2n}^{(N)}\leq \left(\frac{2\alpha}{\nu}+\eps\right)\sum_{i=N}^{n} x_{n-i,n}^{-1}(x_{n-i,n}-x_{n-i-1,n}),
\end{equation}
consequently, if one proves
\begin{equation}\label{eq:equiv-sum-integral}
\sum_{i=N}^{n} x_{n-i,n}^{-1}(x_{n-i,n}-x_{n-i-1,n})\sim \int_{x_{-1,n}}^{x_{n-N,n}} \frac{1}{x}dx,\quad \text{ as }n\to\infty,
\end{equation}
then, by using Lemma \ref{lem:conv-sum-phi}, one obtains
\begin{align}
\sum_{i=N}^{n} x_{n-i,n}^{-1}(x_{n-i,n}-x_{n-i-1,n})&\sim \int_{x_{-1,n}}^{x_{n-N,n}} \frac{1}{x}dx\nonumber\\
&=\log\left(1+\mu_n(1-s_{n+1})\sum_{l=N}^n\frac{\varphi_l(f_{l,n}(s_{n+1}))}{\mu_{l-1}}\right)\nonumber\\
&\sim\log\left(1+\frac{\mu_n\lambda}{2a_{n+1}}\sum_{l=N}^n\frac{\nu_l}{\mu_{l-1}}\right)\nonumber\\
&\to \log(1+\lambda),\label{eq:lim-sum-integ}
\end{align}
as $n\to\infty$. Observe that \eqref{eq:equiv-sum-integral} follows from the integral test for convergence and the fact that given $\delta>0$, there exists $N_2=N_{2,\delta}\in\N$ such that
\begin{equation}\label{eq:cond-equiv-integ}
\sup_{0\leq i\leq n}\left|\frac{x_{n-i,n}}{x_{n-i-1,n}}-1\right|<\delta,\quad \text{ for each }n\geq N_2.
\end{equation}
Indeed, from \eqref{eq:unif-conv-phi} one has that the sequence $\{\varphi_i(f_{i,n}(s_{n+1}))\}_{0\leq i\leq n; n\in\N_0}$ is bounded and
\begin{equation*}
\left|\frac{x_{n-i,n}}{x_{n-i-1,n}}-1\right|=\frac{\varphi_i(f_{i,n}(s_{n+1}))}{\frac{\mu_{i-1}}{\mu_n(1-s_{n+1})}+\mu_{i-1}\sum_{l=i+1}^n\frac{\varphi_l(f_{l,n}(s_{n+1}))}{\mu_{l-1}}}\leq
\frac{K_1\mu_n(1-s_{n+1})}{\mu_{i-1}},
\end{equation*}
where $K_1>0$ is an upper bound of the sequence $\{\varphi_i(f_{i,n}(s_{n+1}))\}_{0\leq i\leq n; n\in\N_0}$. With similar arguments to those used before to establish \eqref{eq:proof-inf-sup}, \eqref{eq:cond-equiv-integ} yields.

Thus, from \eqref{eq:bounds-S2n} and \eqref{eq:lim-sum-integ},
\begin{equation*}
\left(\frac{2\alpha}{\nu}-\varepsilon\right)\log(1+\lambda)\leq \liminf_{n\to\infty} S_{2n}^{(J_0)}\leq \limsup_{n\to\infty} S_{2n}^{(J_0)}\leq \left(\frac{2\alpha}{\nu}+\varepsilon\right)\log(1+\lambda)
\end{equation*}
hence, $\lim_{n\to\infty} S_{2n}^{(J_0)}=\log (1+\lambda)^{\frac{2\alpha}{\nu}}$. For $S_{1n}^{(J_0)}$, one has $S_{1n}^{(J_0)}=o(1)$, as $n\to\infty$; indeed,
\begin{equation*}
S_{1n}^{(J_0)}=\sum_{i=0}^{J_0-1}\alpha_i(1-f_{i-1,n}(s_{n+1}))\leq J_0\max_{0\leq i<J_0}\left\{\frac{\alpha_i}{\mu_{i-1}}\right\}\mu_n(1-s_{n+1})\to 0,
\end{equation*}
which together with the previous limit gives us \eqref{eq:result1-thm-asym-critical}.

Finally, to establish \eqref{eq:result2-thm-asym-critical}, we have that
%\begin{equation}\label{eq:aprox-error-term}
%\sum_{i=0}^n \left(\frac{h_i''(\xi_{in})}{h_i(\xi_{in})}-\frac{h_i'(\xi_{in})^2}{h_i(\xi_{in})^2}\right)(1-f_{i-1,n}(s_{n+1}))^2=o\left(\sum_{i=0}^{n}\alpha_i(1-f_{i-1,n}(s_{n+1}))\right),\quad\text{ as }n\to\infty.
%\end{equation}
%Indeed,
\begin{align*}
\sum_{i=0}^n &\left|\frac{h_i''(\xi_{in})}{h_i(\xi_{in})}-\frac{h_i'(\xi_{in})^2}{h_i(\xi_{in})^2}\right|(1-f_{i-1,n}(s_{n+1}))^2\leq\\
&\leq \left(\sum_{i=0}^n \alpha_i(1-f_{i-1,n}(s_{n+1}))\right)\sup_{0\leq i\leq n}\left\{\frac{1-f_{i-1,n}(s_{n+1})}{\alpha_i}\left|\frac{h_i''(\xi_{in})}{h_i(\xi_{in})}-\frac{h_i'(\xi_{in})^2}{h_i(\xi_{in})^2}\right|\right\},
\end{align*}
and
\begin{align*}
\sup_{0\leq i\leq n}\bigg\{\frac{1-f_{i-1,n}(s_{n+1})}{\alpha_i}&\bigg|\frac{h_i''(\xi_{in})}{h_i(\xi_{in})}-\frac{h_i'(\xi_{in})^2}{h_i(\xi_{in})^2}\bigg|\bigg\}\leq\\ &\leq\bigg(\frac{K_2\sup_{n\in\N_0}h_n''(1)}{\tau}+\frac{K_3}{\tau^2}\bigg)\sup_{0\leq i\leq n} \{1-f_{i,n}(s_{n+1})\},
\end{align*}
with $K_2>0$ and $K_3>0$ being upper bounds of the sequences $\{1/\alpha_{n}\}_{n\in\N_0}$ and $\{\alpha_{n}\}_{n\in\N_0}$, respectively. Notice that $\sup_{n\in\N_0} \beta_n^2<\infty$ implies $\sup_{n\in\N_0}h_n''(1)<\infty$.  Furthermore, from \eqref{eq:lim-sup-proof-asymp}, the convergence $\sup_{0\leq i\leq n} \{1-f_{i-1,n}(s_{n+1})\}\to 0$, as $n\to\infty$, is deduced, and as a consequence, using \eqref{eq:result1-thm-asym-critical},  \eqref{eq:result2-thm-asym-critical} is obtained and the proof is finished.
%To prove this assertion, we apply again a second-order Taylor series expansion
%\begin{align*}
%1-f_{i-1,n}(s_{n+1})&\leq f_{i-1,n}''(1)(1-s_{n+1})+\frac{f_{i-1,n}''(1)}{2}(1-s_{n+1})^2\\
%&\leq \frac{\mu_n}{\mu_{i-1}}(1-s_{n+1})+\frac{\mu_n^2}{2\mu_{i-1}}(1-s_{n+1})^2\sum_{l=i}^n \frac{\nu_l}{\mu_{l-1}}\\
%&\leq 2K_1(1+K_1)\left(\mu_{i-1}\sum_{k=0}^{i-1}\frac{\nu_k}{\mu_{k-1}}\right)^{-1}\to 0,\quad\text{ as }i\to\infty,
%\end{align*}
%and hence, given $\gamma>0$, there exists $I_3=I_{3,\gamma}\in\N$, such that $\sup_{I_3\leq i\leq n} (1-f_{i-1,n}(s_{n+1}))<\gamma/2$, for any $n\geq I_3+1$. On the other hand, since
%\begin{align*}
%\sup_{0\leq i\leq I_3-1} (1-f_{i-1,n}(s_{n+1}))&=\sup_{0\leq i\leq I_3-1} \frac{1}{\frac{\mu_{i-1}}{\mu_n(1-s_{n+1})}+\mu_{i-1}\sum_{l=i}^n\frac{\varphi_l(f_{l,n}(s_{n+1}))}{\mu_{l-1}}}\\
%&\leq \max_{0\leq i\leq I_3-1}\left\{\frac{1}{\mu_{i-1}}\right\}\mu_n(1-s_{n+1})\to 0,\qquad\text{  as }n\to\infty,
%\end{align*}
%there exists $N_5=N_{5,\gamma}\in\N$ such that $\sup_{0\leq i\leq I_3-1} (1-f_{i-1,n}(s_{n+1}))<\gamma/2$, for $n\geq N_5$, and consequently,
%\begin{equation*}
%\sup_{0\leq i\leq n} \{1-f_{i-1,n}(s_{n+1}))\}<\gamma,\quad\text{ for }n\geq\max\{N_5,I_3+1\}.
%\end{equation*}

%Finally, from \eqref{eq:aprox-Fn},  \eqref{eq:lim-Sn} and \eqref{eq:aprox-error-term}, we have
%\begin{equation*}
%F_{n+1}(s_{n+1})=\exp\left\{\log(1+\lambda)^{-\frac{2\alpha}{\nu}}+o\left(\log(1+\lambda)^{-\frac{2\alpha}{\nu}}\right)\right\}\to (1+\lambda)^{-\frac{2\alpha}{\nu}},\quad\text{ as }n\to\infty.
%\end{equation*}
\end{Prf}
\section*{Acknowledgements}

This manuscript was initiated while Carmen Minuesa was a visiting PhD student at the Institute of Mathematics, Goethe University Frankfurt, in Frankfurt am Main, and she is grateful for the hospitality and collaboration.

This research has been supported by the Spanish Ministerio de Educaci\'on, Cultura y Deporte under the Grants FPU13/03213 and EST16/00404, the Spanish Ministerio de Econom\'ia y Competitividad under the Grant MTM2015-70522-P, the Junta de Extremadura under the Grants IB16099 and GR18103 and the Fondo Europeo de Desarrollo Regional.

A revision of this version of the article has been accepted for publication, after peer review and is subject to Taylor \& Francis Ltd.’s terms of use, but is not the Version of Record and does not reflect post-acceptance improvements, or any corrections. The Version of Record is available online at: \url{https://doi.org/10.1080/15326349.2019.1575754}

%%%%%%%%%%%%%%%%%%%%%%%%%%%%%%%%%%%%%%%%%%%%%%%%%%%%%%%%%%%%%%%%%%%%%

%%%%%%%%%%%%%%%%%%%%%%%%%%% REFERENCES %%%%%%%%%%%%%%%%%%%%%%%%%%%%%%%%
%%\addtocontents{toc}{\protect\setcounter{tocdepth}{1}}
%\addcontentsline{toc}{section}{References}
%\bibliographystyle{plain}
%\bibliography{bibtes_tesis}

\begin{thebibliography}{10}

\bibitem{Athreya-Ney}
K.~B. Athreya and P.~E. Ney.
\newblock {\em Branching processes}.
\newblock Springer, 1972.

\bibitem{Souza-1995}
J.~C. D'Souza.
\newblock The extinction time of the inhomogeneous branching process.
\newblock In C.~C. Heyde, editor, {\em Branching Processes. Proceedings of the
  First World Congress}, volume~99 of {\em Lecture Notes in Statistics}, pages
  106--117. Springer New York, 1995.

\bibitem{Fearn-1972}
D.~H. Fearn.
\newblock {Galton-Watson processes with generation dependence}.
\newblock In L.~M. {Le Cam}, J.~Neyman, and E.~L. Scott, editors, {\em
  Proceedings of the Sixth Berkeley Symposium on Mathematical Statistics and
  Probability, Volume 4: Biology and Health}, number~1, pages 159--172.
  University of California Press, 1972.

\bibitem{Foster-Williamson-1971}
J.~H. Foster and J.~A. Williamson.
\newblock {Limit theorems for the Galton-Watson process with time-dependent
  immigration}.
\newblock {\em Zeitschrift f{\"{u}}r Wahrscheinlichkeitstheorie und Verwandte
  Gebiete}, 20(3):227--235, 1971.

\bibitem{Gao-Zhang-2015}
Z.~Gao and Y.~Zhang.
\newblock Limit theorems for a {G}alton--{W}atson process with immigration in
  varying environments.
\newblock {\em Bulletin of Malaysian Mathematical Science Society},
  38:1551--1573, 2015.

\bibitem{CBPs}
M.~Gonz\'alez, I.~del Puerto, and G.~P. Yanev.
\newblock {\em Controlled Branching Processes}.
\newblock ISTE Ltd and John Wiley and Sons, Inc., 2018.

\bibitem{art-Dposterior}
M.~Gonz\'alez, C.~Minuesa, I.~del Puerto, and A.~N. Vidyashankar.
\newblock Robust estimation in controlled branching processes: Bayesian
  estimators via disparities.
\newblock arXiv:1802.05917:1--62, 2018.

\bibitem{LMJ}
M.~Gonz{\'a}lez, C.~Minuesa, M.~Mota, I.~del Puerto, and A.~Ramos.
\newblock An inhomogeneous controlled branching process.
\newblock {\em Lithuanian Mathematical Journal}, 55(1):61--71, 2015.

\bibitem{Grishechkin-1992}
S.~A. Grishechkin.
\newblock {Multiclass batch arrival retrial queues analyzed as branching
  processes with immigration}.
\newblock {\em Queueing Systems}, 11:395--418, 1992.

\bibitem{Haccou}
P.~Haccou, P.~Jagers, and V.~Vatutin.
\newblock {\em Branching processes: variation, growth and extinction of
  populations}.
\newblock Cambridge Studies in Adaptative Dynamics. Cambridge University Press,
  2005.

\bibitem{Heathcote-1965}
C.~R. Heathcote.
\newblock A branching process allowing immigration.
\newblock {\em Journal of the Royal Statistical Society. Series B
  (Methodological)}, 27(1):138--143, 1965.

\bibitem{Heathcote-1966}
C.~R. Heathcote.
\newblock Corrections and comments on the paper ``{A} branching process
  allowing immigration''.
\newblock {\em Journal of the Royal Statistical Society. Series B
  (Methodological)}, 28(1):213--217, 1966.

\bibitem{Hyrien-Yanev-Jordan-2015}
O.~Hyrien, N.~M. Yanev, and C.~T. Jordan.
\newblock A test of homogeneity for age-dependent branching processes with
  immigration.
\newblock {\em Electronic Journal of Statistics}, 9(1):898--925, 2015.

\bibitem{ispany-2016}
M.~Isp\'any.
\newblock Some asymptotic results for strongly critical branching processes
  with immigration in varying environment.
\newblock In I.~del Puerto, M.~Gonz\'{a}lez, C.~Guti\'{e}rrez, R.~Mart\'{i}nez,
  C.~Minuesa, M.~Molina, M.~Mota, and A.~Ramos, editors, {\em Branching
  Processes and Their Applications}, volume 219 of {\em Lecture Notes in
  Statistics}, pages 77--95. Springer, 2016.

\bibitem{jagers-68}
P.~Jagers.
\newblock Age-dependent branching processes allowing immigration.
\newblock {\em Theory of Probability and its Applications}, 13(2):225--236,
  1968.

\bibitem{kersting-2017}
G.~Kersting.
\newblock A unifying approach to branching processes in varying environment.
\newblock {\em arXiv:1703.01960}, pages 1--23, 2017.

\bibitem{BPVEs}
G.~Kersting and V.~Vatutin.
\newblock {\em Discrete Time Branching Processes in Random Environment}.
\newblock ISTE Ltd and John Wiley and Sons, Inc., 2017.

\bibitem{Key-1987}
E.~S. Key.
\newblock {Limiting Distributions and Regeneration Times for Multitype
  Branching Processes with Immigration in a Random Environment}.
\newblock {\em The Annals of Probability}, 15(1):344--353, 1987.

\bibitem{MaCphee-Schuh-1983}
I.~M. MaCphee and H.~J. Schuh.
\newblock A {G}alton-{W}atson branching process in varying environments with
  essentially constant offspring means and two rates of growth.
\newblock {\em Australian Journal of Statistics}, 25(2):329--338, 1983.

\bibitem{Mitov-Omey-2014}
K.~V. Mitov and E.~Omey.
\newblock {A branching process with immigration in varying environments}.
\newblock {\em Communications in Statistics - Theory and Methods},
  43(24):5211--5225, 2014.

\bibitem{Nagaev-Khan-1980}
S.~V. Nagaev and L.~V. Khan.
\newblock {Limit theorems for a critical Galton-Watson branching process with
  migration}.
\newblock {\em Theory of Probability and its Applications}, 25(3):514--525,
  1980.

\bibitem{Pakes-1971}
A.~G. Pakes.
\newblock A branching process with a state dependent immigration component.
\newblock {\em Advances in Applied Probability}, 3(2):301--314, 1971.

\bibitem{R}
{R Core Team}.
\newblock {\em {R}: A language and environment for statistical computing}.
\newblock R Foundation for Statistical Computing, Vienna, Austria, 2017.

\bibitem{Rahimov-1987}
I.~Rahimov.
\newblock Critical branching processes with infinite variance and decreasing
  immigration.
\newblock {\em Theory of Probability and its Applications}, 31(1):88--100,
  1987.

\bibitem{Rahimov-2007}
I.~Rahimov.
\newblock Functional limit theorems for critical processes with immigration.
\newblock {\em Advances in Applied Probability}, 39(4):1054--1069, 2007.

\bibitem{Rahimov-2017}
I.~Rahimov.
\newblock Statistical inference for partially observed branching processes with
  immigration.
\newblock {\em Journal of Applied Probability}, 54(1):82–95, 2017.

\bibitem{Rahimov-2018}
I.~Rahimov.
\newblock Estimation of the mean in partially observed branching processes with
  general immigration.
\newblock {\em Statistical Inference for Stochastic Processes}, 2018.

\bibitem{Sebastyanov-1957}
B.~A. Sevast'yanov.
\newblock Limit theorems for branching random processes of special type.
\newblock {\em Theory of Probability and its Applications}, 2(3):321--331,
  1957.
\newblock In Russian.

\bibitem{Yakovlev-Yanev-2006}
A.~Yakovlev and N.~M. Yanev.
\newblock Branching stochastic processes with immigration in analysis of
  renewing cell populations.
\newblock {\em Mathematical biosciences}, 203(1):37--63, 2006.

\bibitem{Yakovlev-Stoimenova-Yanev-2008}
A.~Y. Yakovlev, V.~K. Stoimenova, and N.~M. Yanev.
\newblock Branching processes as models of progenitor cell populations and
  estimation of the offspring distributions.
\newblock {\em Journal of the American Statistical Association},
  103(484):1357--1366, 2008.

\end{thebibliography}
%%%%%%%%%%%%%%%%%%%%%%%%%%%%%%%%%%%%%%%%%%%%%%%%%%%%%%%%%%%%%%%%%%%%%%

\end{document}